\documentclass[12pt]{article}
\usepackage{amssymb,amsmath}

\newcommand{\Z}{{\mathbb Z}}
\newcommand{\F}{\mathcal{F}}
\newcommand{\G}{\mathcal{G}}

\newcommand{\N}{{\mathbb N}}

\newcommand{\OO}{{\cal O}}
\newcommand{\M}{\mathcal{M}}

\renewcommand{\char}{\mbox{char}}

\newcommand{\id}{\mathrm{id}}
\newcommand{\norm}{\mathrm{N}}

\newcommand{\Kbar}{\overline{K}}
\newcommand{\epsilonbar}{\overline{\epsilon}}
\newcommand{\ubar}{\overline{u}}
\newcommand{\Phibar}{\overline{\Phi}}
\newcommand{\proof}{\noindent{\em Proof: }}
\newcommand{\qed}{\hspace{\fill}$\square$}
\newcommand{\ra}{\rightarrow}

\newcommand{\dst}{\displaystyle}

\newtheorem{theorem}{Theorem}
\newtheorem{lemma}[theorem]{Lemma}
\newtheorem{prop}[theorem]{Proposition}
\newtheorem{cor}[theorem]{Corollary}

\newenvironment{remark}{\noindent\refstepcounter{theorem}{\bf
Remark \thesection.\arabic{theorem}} }{}

\numberwithin{equation}{section}
\numberwithin{theorem}{section}

\title{Indices of inseparability in towers of field
extensions}
\author{Kevin Keating \\
Department of Mathematics \\
University of Florida \\
Gainesville, FL 32611 \\
USA \\[.2cm]
{\tt keating@ufl.edu}}

\begin{document}

\maketitle

\begin{abstract}
\noindent
Let $K$ be a local field whose residue field has
characteristic $p$ and let $L/K$ be a finite separable
totally ramified extension of degree $n=ap^{\nu}$.  The
indices of inseparability $i_0,i_1,\dots,i_{\nu}$ of
$L/K$ were defined by Fried in the case $\char(K)=p$ and
by Heiermann in the case $\char(K)=0$; they give a
refinement of the usual ramification data for $L/K$.
The indices of inseparability can be used to construct
``generalized Hasse-Herbrand functions'' $\phi_{L/K}^j$
for $0\le j\le\nu$.  In this paper we give an
interpretation of the values $\phi_{L/K}^j(c)$ for
nonnegative integers $c$.  We use this interpretation to
study the behavior of generalized Hasse-Herbrand
functions in towers of field extensions.
\end{abstract}

\section{Introduction}

     Let $K$ be a local field whose residue field
$\Kbar$ is a perfect field of characteristic $p$, and
let $K^{sep}$ be a separable closure of $K$.  Let
$L/K$ be a finite totally ramified subextension of
$K^{sep}/K$.  The {\em indices of inseparability} of
$L/K$ were defined by Fried \cite{fried} in the case
$\char(K)=p$, and by Heiermann \cite{heier} in the case
$\char(K)=0$.  The indices of inseparability
of $L/K$ determine the ramification data of $L/K$ (as
defined for instance in Chapter IV of \cite{cl}), but
the ramification data does not always determine the
indices of inseparability.  Therefore the indices of
inseparability of $L/K$ may be viewed as a refinement of
the usual ramification data of $L/K$.

     Let $\pi_K$, $\pi_L$ be uniformizers for $K$, $L$.
The most natural definition of the ramification data
of $L/K$ is based on the valuations of
$\sigma(\pi_L)-\pi_L$ for $K$-embeddings
$\sigma:L\ra K^{sep}$; this is the approach taken in
Serre's book \cite{cl}.  The ramification data can also
be defined in terms of the relation between the norm
map $\norm_{L/K}$ and the filtrations of the unit groups
of $L$ and
$K$, as in Fesenko-Vostokov \cite{FV}.  This approach
can be used to derive the well-known relation between
higher ramification theory and class field theory.
Finally, the ramification data can be computed by
expressing $\pi_K$ as a power series in $\pi_L$ with
coefficients in the set $R$ of Teichm\"uller
representatives for $\Kbar$.  This third approach, which
is used by Fried and Heiermann, makes clear the connection
between ramification data and the indices of
inseparability.

     Heiermann \cite{heier} defined ``generalized
Hasse-Herbrand functions'' $\phi_{L/K}^j$ for
$0\le j\le\nu$.  In Section~\ref{HH} we give an
interpretation of the values
$\phi_{L/K}^j(c)$ of these functions at nonnegative
integers $c$.  This leads to an alternative
definition of the indices of inseparability which is
closely related to the third method for defining the
ramification data.  In Section~\ref{towers} we consider a
tower of finite totally ramified separable extensions
$M/L/K$.  We use our interpretation of the
values $\phi_{L/K}^j(c)$ to study the relations between
the generalized Hasse-Herbrand functions of $L/K$,
$M/L$, and $M/K$. \medskip

\noindent{\bf Notation} \\[\medskipamount]
$\N_0=\N\cup\{0\}=\{0,1,2,\dots\}$ \\
$v_p=p$-adic valuation on $\Z$ \\
$K=$ local field with perfect residue field $\Kbar$ of
characteristic $p>0$ \\
$K^{sep}=$ separable closure of $K$ \\
$v_K=$ valuation on $K^{sep}$ normalized so that
$v_K(K^{\times})=\Z$ \\
$\OO_K=\{\alpha\in K:v_K(\alpha)\ge0\}=$ ring of integers
of $K$ \\
$\pi_K=$ uniformizer for $K$ \\
$\M_K=\pi_K\OO_K=$ maximal ideal of $\OO_K$ \\
$R=$ set of Teichm\"uller representatives for $\Kbar$ \\
$L/K=$ finite totally ramified subextension of
$K^{sep}/K$ of degree $n>1$, with $v_p(n)=\nu$ \\
$M/L=$ finite totally ramified subextension of
$K^{sep}/L$ of degree $m>1$, with $v_p(m)=\mu$ \\[.1cm]
$v_K$, $\OO_K$, $\pi_K$, and $\M_K$ have natural analogs
for $L$ and $M$

\section{Generalized Hasse-Herbrand functions}
\label{HH}

     We begin by recalling the definition of the indices
of inseparability $i_j$ ($0\le j\le\nu$) for a
nontrivial totally ramified separable extension $L/K$ of
degree $n=ap^{\nu}$, as formulated by Heiermann
\cite{heier}.  Let $R\subset\OO_K$ be the set of
Teichm\"uller representatives for $\Kbar$.  Then there
is a unique series
$\dst\hat{\F}(X)=\sum_{h=0}^{\infty}a_hX^{h+n}$
with coefficients in $R$ such that $\pi_K=\hat{\F}(\pi_L)$.
For $0\le j\le\nu$ set
\begin{equation} \label{tildeij}
\tilde{\imath}_j=\min\{h\ge0:v_p(h+n)\le j,\;a_h\not=0\}.
\end{equation}
If $\char(K)=0$ it may happen that $a_h=0$ for all $h\ge0$
such that $v_p(h+n)\le j$, in which case we set
$\tilde{\imath}_j=\infty$.  The indices of inseparability
are defined recursively in terms of $\tilde{\imath}_j$ by
$i_{\nu}=\tilde{\imath}_{\nu}=0$ and
$i_j=\min\{\tilde{\imath}_j,i_{j+1}+v_L(p)\}$ for
$j=\nu-1,\dots,1,0$.  Thus
\begin{equation}
i_j=\min\{\tilde{\imath}_{j_1}+(j_1-j)v_L(p):j\le j_1\le\nu\}.
\end{equation}

     It follows from the definitions that
$0=i_{\nu}<i_{\nu-1}\le i_{\nu-1}\le\dots\le i_0$.  If
$\char(K)=p$ then $v_L(p)=\infty$, so
$i_j=\tilde{\imath}_j$ in this case.  If $\char(K)=0$
then $\tilde{\imath}_j$ can depend on the choice of
$\pi_L$, and it is not obvious that $i_j$ is a
well-defined invariant of the extension $L/K$.  We will
have more to say about this issue in Remark~\ref{well}.

     Following \cite[(4.4)]{heier}, for $0\le j\le\nu$
we define functions
$\tilde{\phi}_{L/K}^j:[0,\infty)\ra[0,\infty)$ by
$\tilde{\phi}_{L/K}^j(x)=i_j+p^jx$.  The generalized
Hasse-Herbrand functions
$\phi_{L/K}^j:[0,\infty)\ra[0,\infty)$ are then defined
by
\begin{equation} \label{phik}
\phi_{L/K}^j(x)
=\min\{\tilde{\phi}_{L/K}^{j_0}(x):0\le j_0\le j\}.
\end{equation}
Hence we have 
$\phi_{L/K}^j(x)\le\phi_{L/K}^{j'}(x)$ for
$0\le j'\le j$.
Let $\phi_{L/K}:[0,\infty)\ra[0,\infty)$ be the usual
Hasse-Herbrand function, as defined for instance in
Chapter IV of \cite{cl}.  Then by
\cite[Cor.\,6.11]{heier} we have
$\phi_{L/K}^{\nu}(x)=n\phi_{L/K}(x)$.

     In order to reformulate the definition of
$\phi_{L/K}^j(x)$ we will use the following elementary
fact about binomial coefficients, which is proved in
\cite[Lemma~5.6]{heier}.

\begin{lemma} \label{val}
Let $b\ge c\ge1$.  Then
$\dst v_p\left(\binom{b}{c}\right)\ge v_p(b)-v_p(c)$,
with equality if ${v_p(b)\ge v_p(c)}$ and $c$ is a power
of $p$.
\end{lemma}

\begin{prop} \label{newdef}
For $0\le j\le\nu$ and $x\ge0$ we have
\begin{equation} \nonumber
\phi_{L/K}^j(x)=
\min\left\{h+v_L\left(\binom{h+n}{p^{j_0}}\right)+p^{j_0}x:
0\le j_0\le j,\;a_h\not=0\right\}.
\end{equation}
\end{prop}

\proof Using (\ref{tildeij})--(\ref{phik}) we get
\begin{align*}
\phi_{L/K}^j(x)&=\min\{h+(j_1-j_0)v_L(p)+p^{j_0}x: \\
&\hspace*{5cm}0\le j_0\le j,\;j_0\le j_1\le\nu,\;
v_p(h+n)\le j_1,\;a_h\not=0\}.
\end{align*}
If $j_0>v_p(h+n)$ then we can replace $j_0$ with $j_0-1$
and $j_1$ with $j_1-1$ without increasing the value of
$h+(j_1-j_0)v_L(p)+p^{j_0}x$.  Hence we may assume
$j_0\le v_p(h+n)$ and $j_1=v_p(h+n)$.  It follows that
\begin{align*}
\phi_{L/K}^j(x)&=\min\{h+(v_p(h+n)-j_0)v_L(p)+p^{j_0}x: 
0\le j_0\le j,\;j_0\le v_p(h+n),\;a_h\not=0\} \\
&=\min\left\{h+v_L\left(\binom{h+n}{p^{j_0}}\right)+p^{j_0}x:
0\le j_0\le j,\,j_0\le v_p(h+n),\;a_h\not=0\right\} \\
&=\min\left\{h+v_L\left(\binom{h+n}{p^{j_0}}\right)+p^{j_0}x:
0\le j_0\le j,\;a_h\not=0\right\},
\end{align*}
where the second and third equalities follow from
Lemma~\ref{val}. \qed \medskip

     For $d\ge0$ set $B_d=\OO_L/\M_L^{n+d}$ and let
$A_d=(\OO_K+\M_L^{n+d})/\M_L^{n+d}$ be the image of
$\OO_K$ in $B_d$.  For $0\le j\le\nu$ set
$B_d[\epsilon_j]=B_d[\epsilon]/(\epsilon^{p^{j+1}})$, so
that $\epsilon_j=\epsilon+(\epsilon^{p^{j+1}})$ satisfies
$\epsilon_j^{p^{j+1}}=0$.

\begin{prop} \label{series}
Let $0\le j\le\nu$, let $d\ge c\ge0$, and let $u\in
\OO_L[\epsilon_j]^{\times}$.  Choose
$F(X)\in X^n\cdot\OO_K[[X]]$ such that
$F(\pi_L)=\pi_K$.  Then the following are equivalent:
\begin{enumerate}
\item $F(\pi_L+u\pi_L^{c+1}\epsilon_j)\equiv\pi_K
\pmod{\pi_L^{n+d}}$.
\item There exists an $A_d$-algebra homomorphism
$s_d:B_d\ra B_d[\epsilon_j]$ such that
$s_d(\pi_L)=\pi_L+u\pi_L^{c+1}\epsilon_j$.
\item There exists an $A_d$-algebra homomorphism
$s_d:B_d\ra B_d[\epsilon_j]$ such that
\begin{align*}
s_d&\equiv\id_{B_d}\pmod{\pi_L^{c+1}\epsilon_j} \\
s_d&\not\equiv\id_{B_d}
\pmod{\pi_L^{c+1}\epsilon_j\cdot(\pi_L,\epsilon_j)}.
\end{align*}
\end{enumerate}
\end{prop}

\proof Suppose Condition~1 holds.  Let
$\tilde{u}(X,\epsilon_j)$ be an element of
$\OO_K[[X]][\epsilon_j]$ such that
$\tilde{u}(\pi_L,\epsilon_j)=u$.  Since $F(0)=0$ the
Weierstrass polynomial of $F(X)-\pi_K$ is the minimum
polynomial of $\pi_L$ over $K$.  Therefore
$\OO_L\cong\OO_K[[X]]/(F(X)-\pi_K)$.  It follows that
the $\OO_K$-algebra homomorphism
$\tilde{s}:\OO_K[[X]]\ra\OO_K[[X]][\epsilon_j]$ defined
by $\tilde{s}(X)=X+\tilde{u}X^{c+1}\epsilon_j$ induces an
$A_d$-algebra homomorphism $s_d:B_d\ra B_d[\epsilon_j]$
such that $s_d(\pi_L)=\pi_L+u\pi_L^{c+1}\epsilon_j$.
Therefore Condition~2 holds.  On the other hand, if
Condition~2 holds then applying the homomorphism $s_d$
to the congruence
$F(\pi_L)\equiv\pi_K\pmod{\pi_L^{n+d}}$ gives
Condition~1.  Hence the first two conditions are equivalent.
Suppose Condition~2 holds.  Since $d\ge c$ and $n\ge2$
we see that $s_d$ satisfies the requirements of
Condition~3.  Suppose Condition~3 holds.  Then
$s_d(\pi_L)=\pi_L+v\pi_L^{c+1}\epsilon_j$ for some
$v\in B_d[\epsilon_j]^{\times}$.  Let
$\gamma:B_d[\epsilon_j]\ra B_d[\epsilon_j]$ be the
$B_d$-algebra homomorphism such that
$\gamma(\epsilon_j)=uv^{-1}\epsilon_j$, and define
$s_d':B_d\ra B_d[\epsilon_j]$ by $s_d'=\gamma\circ s_d$.
Then $s_d'$ satisfies the requirements of
Condition~2.~\qed \medskip

     The assumptions on $F(X)$ imply that
$F(\pi_L+u\pi_L^{c+1}\epsilon_j)\equiv\pi_K
\pmod{\pi_L^{n+c}}$.  Therefore the conditions of the
proposition are satisfied when $d=c$.  On the other
hand, since $L/K$ is separable we have
$F(\pi_L+u\pi_L^{c+1}\epsilon_j)\not=\pi_K$.  Hence for
$d$ sufficiently large the conditions in the proposition
are not satisfied.
We define a function $\Phi_{L/K}^j:\N_0\ra\N_0$ by
setting $\Phi_{L/K}^j(c)$ equal to the largest integer
$d$ satisfying the equivalent conditions of
Proposition~\ref{series}.  By Condition~3 we see that
this definition does not depend on the choice of
$\pi_L$, $u$, or $F$.

     We now show that $\Phi_{L/K}^j$ and $\phi_{L/K}^j$
agree on nonnegative integers.  This gives an
alternative description of the restriction of
$\phi_{L/K}^j$ to $\N_0$ which does not depend on the
indices of inseparability.

\begin{prop} \label{defn}
For $c\in\N_0$ we have $\Phi_{L/K}^j(c)=\phi_{L/K}^j(c)$.
\end{prop}

\proof Let $c\in\N_0$.  Since $\hat{\F}(X)$ satisfies
the hypotheses for $F(X)$ in Proposition~\ref{series},
$\Phi_{L/K}^j(c)$ is equal to the largest $d\in\N_0$
such that
\begin{equation} \label{Feps}
\hat{\F}(\pi_L+\pi_L^{c+1}\epsilon_j)
\equiv\hat{\F}(\pi_L)\pmod{\pi_L^{n+d}}.
\end{equation}
For $m\ge0$ define
\begin{equation} \nonumber
(D^m\hat{\F})(X)=\sum_{h=0}^{\infty}\,
\binom{h+n}{m}a_hX^{h+n-m}.
\end{equation}
Then
\begin{equation} \nonumber
\hat{\F}(X+\epsilon_j X^{c+1})
=\sum_{m=0}^{p^{j+1}-1}\,(D^m\hat{\F})(X)\cdot
(\epsilon_j X^{c+1})^m.
\end{equation}
Since
$\epsilon_j,\epsilon_j^2,\dots,\epsilon_j^{p^{j+1}-1}$
are linearly independent over $\OO_L$, (\ref{Feps}) holds
if and only if
\begin{equation} \label{Dm}
(D^m\hat{\F})(\pi_L)\cdot\pi_L^{(c+1)m}\in\M_L^{n+d}
\mbox{ for }1\le m<p^{j+1}.
\end{equation}
Hence by Proposition~\ref{newdef} it is sufficient
to prove that (\ref{Dm}) is equivalent to the following:
\begin{equation} \label{pi}
h+v_L\left(\binom{h+n}{p^{j_0}}\right)+cp^{j_0}\ge d
\mbox{ for all $j_0,h$ such that $0\le j_0\le j$ and
$a_h\not=0$.}
\end{equation}

     Assume first that (\ref{pi}) holds.  Choose $m$
such that $1\le m<p^{j+1}$ and write $m=rp^{j_0}$ with
$p\nmid r$ and $j_0\le j$.  Choose $h\ge0$ such that
$a_h\not=0$ and set $l=v_p(h+n)$.  If $m>h+n$ then
$\dst\binom{h+n}{m}=0$, so we have
\begin{equation} \label{MLnd}
\binom{h+n}{m}a_h\pi_L^{h+n-m}\cdot\pi_L^{(c+1)m}
\in\M_L^{n+d}.
\end{equation}
Suppose $m\le h+n$ and $l\ge j_0$.  Using
Lemma~\ref{val} we get
\begin{equation} \nonumber
v_p\left(\binom{h+n}{m}\right)\ge l-j_0
=v_p\left(\binom{h+n}{p^{j_0}}\right).
\end{equation}
Combining this with (\ref{pi}) we get
\begin{equation} \nonumber
h+v_L\left(\binom{h+n}{m}\right)+cm+n\ge
h+v_L\left(\binom{h+n}{p^{j_0}}\right)+cp^{j_0}+n
\ge n+d.
\end{equation}
Hence (\ref{MLnd}) holds in this case.  Finally, suppose
$m\le h+n$ and $l<j_0\le j$.  It follows from
Lemma~\ref{val} that
$\dst v_L\left(\binom{h+n}{p^l}\right)=0$, so by
(\ref{pi}) we have $h+cp^l\ge d$.  Since
$m\ge p^{j_0}>p^l$ we get
\begin{equation} \nonumber
h+v_L\left(\binom{h+n}{m}\right)+cm+n\ge h+cp^l+n\ge n+d.
\end{equation}
Therefore (\ref{MLnd}) holds in this case as well.
It follows that every term in $(D^m\hat{\F})(\pi_L)$
lies in $\M_L^{n+d}$, so (\ref{Dm}) holds.

     Assume conversely that (\ref{Dm}) holds.  Among all
the nonzero terms that occur in any of the series
\begin{equation} \nonumber
(D^{p^i}\hat{\F})(\pi_L)\cdot\pi_L^{(c+1)p^i}
=\sum_{h=0}^{\infty}\,a_h\binom{h+n}{p^i}\pi_L^{h+n+cp^i}
\end{equation}
for $0\le i\le j$ let
$a_{h}\dst\binom{h+n}{p^i}\pi_L^{n+{h}+cp^i}$ be a term
whose $L$-valuation $w$ is minimum.  If $\char(K)=p$
then for each
$m\ge1$ the nonzero terms of $(D^m\hat{\F})(\pi_L)$ have
distinct $L$-valuations, so it follows from (\ref{Dm})
that $w\ge n+d$.  Suppose $\char(K)=0$ and set
$l=v_p(h+n)$.  If $i>l$ then since
$\dst v_L\left(\binom{h+n}{p^l}\right)=0$ we have
\begin{equation} \nonumber
v_L\left(\binom{h+n}{p^l}\pi_L^{n+h+cp^l}\right)\le
v_L\left(\binom{h+n}{p^i}\pi_L^{n+h+cp^i}\right)=w.
\end{equation}
Therefore we may assume $i\le l$.  Since
$\dst v_p\left(\binom{n}{p^i}\right)=\nu-i$ and
$a_0\not=0$ we have $l\le\nu$.  Suppose $w<n+d$.  Then
it follows from (\ref{Dm}) that there
is $h'\not=h$ such that $a_{h'}\not=0$ and
\begin{equation} \label{vL}
v_L\left(\binom{h'+n}{p^i}\pi_L^{n+h'+cp^i}\right)
=v_L\left(\binom{h+n}{p^i}\pi_L^{n+h+cp^i}\right).
\end{equation}
Since $n\mid v_L(p)$ this implies $h'\equiv h\pmod{n}$.
Since $v_p(h+n)\le\nu$ and $v_p(h'+n)\le\nu$ we get
$v_p(h'+n)=v_p(h+n)=l$.  Therefore by Lemma~\ref{val}
we have
\begin{equation} \nonumber
v_p\left(\binom{h'+n}{p^i}\right)
=v_p\left(\binom{h+n}{p^i}\right)=l-i.
\end{equation}
Combining this with (\ref{vL}) gives $h'=h$, a
contradiction.  Therefore $w\ge n+d$ holds in general.
Hence by the minimality of $w$ we get (\ref{pi}). \qed
\medskip

\begin{remark} \label{well}
If $\char(K)=0$ then the value of $\tilde{\imath}_j$
may depend on the choice of uniformizer $\pi_L$ for $L$.
It was proved in \cite[Th.\,7.1]{heier} that $i_j$ is a
well-defined invariant of the extension $L/K$.  This can
also be deduced from Proposition~\ref{defn} by setting
$c=0$.
\end{remark} \medskip

\begin{remark}
Let $0\le j\le\nu$.  Even though the function
$\phi_{L/K}^j:[0,\infty)\ra[0,\infty)$ may not be
determined by its restriction to $\N_0$, it is
determined by the sequence $(i_0,i_1,\dots,i_j)$.
Since $i_{j_0}=\phi_{L/K}^{j_0}(0)$ this implies that
the collection consisting of the restrictions of
$\phi_{L/K}^{j_0}$ to $\N_0$ for $0\le j_0\le j$
determines $\phi_{L/K}^j$.
\end{remark} \medskip

     For $0\le j\le\nu$ let
$B_d[\epsilonbar_j]=B_d[\epsilon]/(\epsilon^{p^j+1})$,
so that $\epsilonbar_j=\epsilon+(\epsilon^{p^j+1})$
satisfies $\epsilonbar_j^{p^j+1}=0$.  Define
$\Phibar_{L/K}^j:\N_0\ra\N_0$ analogously to
$\Phi_{L/K}^j$, using $\epsilonbar_j$ in place of
$\epsilon_j$.  Then the arguments in this section remain
valid with $\epsilon_j,\Phi_{L/K}^j$ replaced by
$\epsilonbar_j,\Phibar_{L/K}^j$.  (In particular, note
that the proof that (\ref{Dm}) implies (\ref{pi}) only
uses the fact that (\ref{Dm}) holds with $m=p^i$ for
$0\le i\le j$.) Hence by Propositions~\ref{series} and
\ref{defn} and their analogs for
$\epsilonbar_j,\Phibar_{L/K}^j$ we get the following:

\begin{cor} \label{ged}
Let $c,d\in\N_0$, let $u\in\OO_L[\epsilon_j]^{\times}$,
and let $\ubar\in\OO_L[\epsilonbar_j]^{\times}$.  Choose
$F(X)\in X^n\cdot\OO_K[[X]]$ such that $F(\pi_L)=\pi_K$.
Then the following are equivalent:
\begin{enumerate}
\item $\phi_{L/K}^j(c)\ge d$,
\item $F(\pi_L+u\pi_L^{c+1}\epsilon_j)\equiv F(\pi_L)
\pmod{\pi_L^{n+d}}$,
\item $F(\pi_L+\ubar\pi_L^{c+1}\epsilonbar_j)\equiv
F(\pi_L)\pmod{\pi_L^{n+d}}$.
\end{enumerate}
\end{cor}

     Some of the proofs in Section~\ref{towers} depend
on ``tame shifts'':

\begin{lemma} \label{unifs}
Let $\pi_L$ be a uniformizer for $L$ and choose a
uniformizer $\pi_K$ for $K$ such that
$\pi_K\equiv\pi_L^n\pmod{\pi_L^{n+1}}$.  Let
$e\ge1$ be relatively prime to $p[L:K]=pn$ and let
$\pi_{K_e}\in K^{sep}$ be a root of $X^e-\pi_K$.  Set
$K_e=K(\pi_{K_e})$ and $L_e=LK_e$.  Then
\begin{enumerate}
\item $K_e/K$ and $L_e/L$ are totally ramified
extensions of degree $e$.
\item There is a uniformizer $\pi_{L_e}$ for $L_e$ such
that $\pi_{L_e}^e=\pi_L$ and
$\pi_{L_e}^n\equiv\pi_{K_e}\pmod{\pi_{L_e}^{n+1}}$.
\item Let $F(X)\in X^n\cdot\OO_K[[X]]$ be such that
$F(\pi_L)=\pi_K$.  Then we can define a series
$F_e(X)=F(X^e)^{1/e}$ with coefficients in $\OO_K$ such
that $F_e(\pi_{L_e})=\pi_{K_e}$.
\end{enumerate}
\end{lemma}

\proof Statement 1 is clear.  Since $e$ and $n$ are
relatively prime there are $s,t\in\Z$ such that
$es+nt=1$.  Then $\tilde{\pi}_{L_e}=\pi_L^s\pi_{K_e}^t$
is a uniformizer for $L_e$ with
$\tilde{\pi}_{L_e}^e\equiv\pi_L
\pmod{\tilde{\pi}_{L_e}^{e+1}}$ and
$\tilde{\pi}_{L_e}^n\equiv\pi_{K_e}
\pmod{\tilde{\pi}_{L_e}^{n+1}}$.
Hence there is a 1-unit $v\in\OO_{L_e}^{\times}$ such
that $\pi_{L_e}=v\tilde{\pi}_{L_e}$ satisfies the
requirements of Statement 2.  To prove Statement 3 we
note that since
$\pi_{L_e}^n\equiv\pi_{K_e}\pmod{\pi_{L_e}^{n+1}}$,
the coefficient $a_0$ in the series
$F(X)=a_0X^n+a_1X^{n+1}+\dots$ is a 1-unit.  Therefore
we may define
\[F_e(X)=F(X^e)^{1/e}
=(a_0X^{ne}+a_1X^{ne+e}+\dots)^{1/e}
=a_0^{1/e}X^n(1+a_0^{-1}a_1X^e+\dots)^{1/e},\]
where $a_0^{1/e}$ is the unique 1-unit in $\OO_K$ whose
$e$th power is $a_0$.  Since
$\pi_{K_e}^e=F(\pi_{L_e}^e)$ and
$\pi_{L_e}^n\equiv\pi_{K_e}\pmod{\pi_{L_e}^{n+1}}$ we
get $F_e(\pi_{L_e})=\pi_{K_e}$.
\qed

\begin{lemma} \label{tame}
Let $K_e$, $L_e$ be as in Lemma~\ref{unifs}.  Then for
$x\ge0$ and $0\le j\le\nu$ we have
\begin{align*}
\tilde{\phi}_{L_e/K_e}^j(x)&=e\tilde{\phi}_{L/K}^j(x/e)
\\
\phi_{L_e/K_e}^j(x)&=e\phi_{L/K}^j(x/e).
\end{align*}
\end{lemma}

\proof It suffices to show that
$ei_0,ei_1,\dots,ei_{\nu}$ are the indices of
inseparability of $L_e/K_e$.  By Proposition~\ref{defn}
this is equivalent to showing that
$\Phi_{L_e/K_e}^j(0)=e\Phi_{L/K}^j(0)$.  Let
$\pi_K$, $\pi_L$, $\pi_{K_e}$, $\pi_{L_e}$,
$F(X)$, $F_e(X)$ satisfy the conditions of
Lemma~\ref{unifs}.  If $\Phi_{L/K}^j(0)\ge d$ then
\begin{alignat*}{2}
F_e(\pi_{L_e}+\pi_{L_e}\epsilon_j)^e
&=F(\pi_L(1+\epsilon_j)^e) \\
&\equiv F(\pi_L)&&\pmod{\pi_L^{n+d}} \\
&\equiv F_e(\pi_{L_e})^e&&\pmod{\pi_L^{n+d}}.
\end{alignat*}
Since $F_e(X)=a_0^{1/e}X^n+\dots$ with $a_0^{1/e}$ a
1-unit, it follows that
\[F_e(\pi_{L_e}+\pi_{L_e}\epsilon_j)
\equiv F_e(\pi_{L_e})\pmod{\pi_{L_e}^{n+de}}.\]
Therefore $\Phi_{L_e/K_e}^j(0)\ge de$.  Conversely,
if $\Phi_{L_e/K_e}^j(0)\ge d$ then
\begin{alignat*}{2}
F(\pi_L+\pi_L\epsilon_j)
&=F_e(\pi_{L_e}(1+\epsilon_j)^{1/e})^e \\
&\equiv F_e(\pi_{L_e})^e
&&\pmod{\pi_K\cdot\pi_{L_e}^d} \\
&\equiv F(\pi_L)&&\pmod{\pi_L^{n+\lceil d/e\rceil}},
\end{alignat*}
and hence $\Phi_{L/K}^j(0)\ge\lceil d/e\rceil$.  By
combining these results we get
$\Phi_{L_e/K_e}^j(0)=e\Phi_{L/K}^j(0)$.~\qed

\section{Towers of extensions} \label{towers}

     In this section we consider a tower $M/L/K$ of
finite totally ramified subextensions of $K^{sep}/K$.
Our goal is to determine relations between the
generalized Hasse-Herbrand functions $\phi_{M/K}^l$ of
the extension $M/K$ and the corresponding functions for
$L/K$ and $M/L$.  It is well-known that the indices of
inseparability of $L/K$ and $M/L$ do not always
determine the indices of inseparability of $M/K$ (see
for instance Example~5.8 in \cite{fm} or Remark~7.8 in
\cite{heier}).  Therefore we cannot expect to obtain a
general formula which expresses $\phi_{M/K}^l$ in terms
of $\phi_{L/K}^j$ and $\phi_{M/L}^k$.  However, we do
get a lower bound for $\phi_{M/K}^l(x)$, and we are
able to show that this lower bound is equal to
$\phi_{M/K}^l(x)$ in certain cases.

     Set $[L:K]=n$, $[M:L]=m$, $\nu=v_p(n)$, and
$\mu=v_p(m)$.  Let $\pi_K$, $\pi_L$, $\pi_M$ be
uniformizers for $K$, $L$, $M$.  Choose
$F(X)\in X^n\cdot\OO_K[[X]]$ such that $F(\pi_L)=\pi_K$
and define
\begin{equation} \nonumber
F^*(\epsilon)=\pi_K^{-1}(F(\pi_L+\pi_L\epsilon)-\pi_K).
\end{equation}
Then $F^*(\epsilon)\in\OO_L[[\epsilon]]$ is uniquely
determined by $L/K$ up to multiplication by an element
of $\OO_L[[\epsilon]]^{\times}$.

     Write
$F^*(\epsilon)=c_1\epsilon+c_2\epsilon^2+\cdots$ and
define the ``valuation function'' of $F^*$ with respect
to $v_K$ by
\begin{equation} \label{Psi}
\Psi_{F^*(\epsilon)}^K(x)=\min\{v_K(c_i)+ix:i\ge1\}
\end{equation}
for $x\in[0,\infty)$.  The graph of
$\Psi_{F^*(\epsilon)}^K$ is the Newton copolygon of
$F^*(\epsilon)$ with respect to $v_K$.  Gross
\cite[Lemma~1.5]{plie} attributes the following
observation to Tate:

\begin{prop} \label{copoly}
For $x\ge0$ we have $\phi_{L/K}(x)=\Psi_{F^*(\epsilon)}^K(x)$.
\end{prop}

     Suppose we also have $G(X)\in X^m\cdot\OO_K[[X]]$
such that $G(\pi_M)=\pi_L$.  Set $H(X)=F(G(X))$.  Then
$H(X)\in X^{nm}\cdot\OO_K[[X]]$ satisfies
$H(\pi_M)=\pi_K$.  It follows that we can use the series
\begin{align*}
G^*(\epsilon)&=\pi_L^{-1}(G(\pi_M+\pi_M\epsilon)-\pi_L) \\
H^*(\epsilon)&=\pi_K^{-1}(H(\pi_M+\pi_M\epsilon)-\pi_K)
\end{align*}
to compute the Hasse-Herbrand functions for the
extensions $M/L$ and $M/K$.  As Lubin points out in
\cite[Th.\,1.6]{eam}, by applying
Proposition~\ref{copoly} to the relation
$H^*(\epsilon)=F^*(G^*(\epsilon))$, we
obtain the well-known composition formula
$\phi_{M/K}=\phi_{L/K}\circ\phi_{M/L}$.

     We wish to extend the results above to apply to the
generalized Hasse-Herbrand functions $\phi_{L/K}^j$.
For $0\le j\le \nu$ let $F^*(\epsilon_j)$ denote the
image of $F^*(\epsilon)$ in
$\OO_L[[\epsilon]]/(\epsilon^{p^{j+1}})\cong\OO_L[\epsilon_j]$.
Alternatively, we may view $F^*(\epsilon_j)$ as the
polynomial obtained by discarding all the terms of
$F^*(\epsilon)$ of degree $\ge p^{j+1}$.  Therefore it
makes sense to consider the valuation function
$\Psi_{F^*(\epsilon_j)}^L(x)$ of $F^*(\epsilon_j)$.

\begin{prop}
$\phi_{L/K}^j(x)=\Psi_{F^*(\epsilon_j)}^L(x)$ for all
$x\in[0,\infty)$.
\end{prop}

\proof We first prove that $\phi_{L/K}^j$ and
$\Psi_{F^*(\epsilon_j)}^L$ agree on $\N_0$.  Let
$d\ge b\ge0$.  Then $\Phi_{L/K}^j(b)\ge d$ if and only
if $F^*(\pi_L^b\epsilon_j)\equiv0\pmod{\pi_L^d}$.  By
(\ref{Psi}) this is equivalent to
$\Psi_{F^*(\epsilon_j)}^L(b)\ge d$.  Since
$\Phi_{L/K}^j$ and $\Psi_{F^*(\epsilon_j)}^L$ map $\N_0$
to $\N_0$, this implies
$\Phi_{L/K}^j(c)=\Psi_{F^*(\epsilon_j)}^L(c)$ for all
$c\in\N_0$.  Using Proposition~\ref{defn} we deduce that
$\phi_{L/K}^j(c)=\Psi_{F^*(\epsilon_j)}^L(c)$ for
$c\in\N_0$.

     Now choose $e\ge1$ relatively prime to
$p[L:K]=pn$.  Let $K_e$, $L_e$, $\pi_K$, $\pi_{K_e}$,
$\pi_L$, $\pi_{L_e}$ satisfy the conditions of
Lemma~\ref{unifs}, and choose $F(X)\in X^n\cdot\OO_K[[X]]$
such that $F(\pi_L)=\pi_K$.  Then $F_e(X)=F(X^e)^{1/e}$
satisfies $F_e(\pi_{L_e})=\pi_{K_e}$.  Let
\begin{align*}
F_e^*(\epsilon)&=\pi_{K_e}^{-1}(F_e(\pi_{L_e}
+\pi_{L_e}\epsilon)-\pi_{K_e}) \\
&=(1+F^*((1+\epsilon)^e-1))^{1/e}-1.
\end{align*}
Then $F_e^*(\epsilon)=\eta^{-1}(F^*(\eta(\epsilon)))$,
where $\eta(\epsilon)=(1+\epsilon)^e-1$ and
$\eta^{-1}(\epsilon)=(1+\epsilon)^{1/e}-1$ have
coefficients in $\OO_K$.  It follows that for
$0\le j\le\nu$ we have
$F_e^*(\epsilon_j)=\eta^{-1}(F^*(\eta(\epsilon_j)))$,
so for $c\in\N_0$ we get
\begin{equation} \nonumber
\Psi_{F_e^*(\epsilon_j)}^{L_e}(c)
=\Psi_{F^*(\epsilon_j)}^{L_e}(c)
=e\Psi_{F^*(\epsilon_j)}^L(c/e).
\end{equation}
By Lemma~\ref{tame} we have
$\phi_{L/K}^j(c/e)=e^{-1}\phi_{L_e/K_e}^j(c)$.  Since
the proposition holds for the extension $L_e/K_e$ with
$x=c$ this implies
\begin{align*}
\phi_{L/K}^j(c/e)&=e^{-1}\Psi_{F_e^*(\epsilon_j)}^{L_e}(c)
=\Psi_{F^*(\epsilon_j)}^L(c/e).
\end{align*}
Since the set $\{c/e:c,e\in\N,\;\gcd(e,pn)=1\}$ is dense
in $[0,\infty)$, and $\phi_{L/K}^j$,
$\Psi_{F^*(\epsilon_j)}^L$ are continuous on
$[0,\infty)$, we conclude that
$\phi_{L/K}^j(x)=\Psi_{F^*(\epsilon_j)}^L(x)$ for all
$x\in[0,\infty)$. \qed \medskip

     Following \cite[(4.4)]{heier}, for $0\le j\le \nu$
and $m\in\N$ we define functions on $[0,\infty)$ by
\begin{align*}
\tilde{\phi}_{L/K}^{j,m}(x)
&=m\tilde{\phi}_{L/K}^j(x/m)=mi_j+p^jx \\
\phi_{L/K}^{j,m}(x)&=m\phi_{L/K}^j(x/m)
=\min\{\tilde{\phi}_{L/K}^{j_0,m}(x):0\le j_0\le j\}.
\end{align*}
For $0\le l\le \nu+\mu$ let
\begin{align*}
\Omega_l&=\{(j,k):0\le j\le \nu,\;0\le k\le \mu,\;j+k=l\},
\end{align*}
and for $x\ge0$ define
\begin{align*}
\lambda_{M/K}^l(x)&=\min\{\phi_{L/K}^{j,m}(\phi_{M/L}^k(x)):
(j,k)\in \Omega_l\} \\
&=\min\{\tilde{\phi}_{L/K}^{j,m}(\tilde{\phi}_{M/L}^k(x)):
(j,k)\in \Omega_{l_0}\text{ for some }0\le l_0\le l\}.
\end{align*}
For $0\le a\le l$ set
\begin{align*}
S_l^a(x)&=\{(j,k)\in \Omega_a:
\tilde{\phi}_{L/K}^{j,m}(\tilde{\phi}_{M/L}^k(x))
=\lambda_{M/K}^l(x)\}.
\end{align*}

\begin{theorem} \label{ge}
Let $0\le l\le\nu+\mu$ and $x\in[0,\infty)$.  Then
\\[\smallskipamount]
(a) $\phi_{M/K}^l(x)\ge\lambda_{M/K}^l(x)$.
\\[\smallskipamount]
(b) Suppose there exists $l_0\le l$ such that
$|S_l^{l_0}(x)|=1$.  Then
$\phi_{M/K}^l(x)=\lambda_{M/K}^l(x)$.
\end{theorem}

     The rest of the paper is devoted to proving this
theorem.  We first consider the cases where $x=c\in\N_0$.
The proof in these cases is based on
Proposition~\ref{defn}.  To get information about
$\Phi_{M/K}^l(c)$ we compute the most significant terms
of $\hat{\F}(\hat{\G}(\pi_M+\pi_M^{c+1}\epsilon))$.

     It follows from Proposition~\ref{defn} that for
$0\le j\le\nu$ we have
\begin{alignat*}{2}
\hat{\F}(\pi_L(1+\epsilon))&\equiv\pi_K
&&\pmod{(\pi_L^{n+i_j},\epsilon^{p^{j+1}})}.
\end{alignat*}
In addition, since $X^n$ divides $\hat{\F}(X)$ we have
\begin{alignat}{2} \label{X}
\hat{\F}(\pi_L(1+\epsilon))&\equiv\pi_K
&&\pmod{\pi_L^n\epsilon}.
\end{alignat}
Hence
\begin{alignat}{2} \label{xeps}
\hat{\F}(\pi_L(1+\epsilon))&\equiv\pi_K
&&\pmod{\pi_L^n\cdot(\pi_L^{i_j},\epsilon^{p^{j+1}})}.
\end{alignat}
Define an ideal in $\OO_L[[\epsilon]]$ by
\begin{align*}
I_{\F}&=(\pi_L^{i_0},\epsilon^{p^1})\cap
(\pi_L^{i_1},\epsilon^{p^2})\cap\dots\cap
(\pi_L^{i_{\nu}},\epsilon^{p^{\nu+1}})\cap
(\epsilon) \\
&=(\pi_L^{i_0}\epsilon^{p^0},\pi_L^{i_1}\epsilon^{p^1},
\dots,\pi_L^{i_{\nu}}\epsilon^{p^{\nu}}).
\end{align*}
It follows from (\ref{X}) and (\ref{xeps}) that
\begin{alignat}{2} \label{Fcong}
\hat{\F}(\pi_L(1+\epsilon))&\equiv\pi_K
&&\pmod{\pi_L^n\cdot I_{\F}}.
\end{alignat}

     Let $i_0',i_1',\dots,i_{\mu}'$ be the indices of
inseparability of $M/L$.  As above we find that
\begin{alignat*}{2}
\hat{\G}(\pi_M(1+\epsilon))&\equiv\pi_L
&&\pmod{\pi_M^m\cdot I_{\G}},
\end{alignat*}
where $I_{\G}$ is the ideal in $\OO_M[[\epsilon]]$
defined by
\begin{align*}
I_{\G}&=
(\pi_M^{i_0'}\epsilon^{p^0},\pi_M^{i_1'}\epsilon^{p^1},\dots,
\pi_M^{i_{\mu}'}\epsilon^{p^{\mu}}).
\end{align*}
By replacing $\epsilon$ with $\pi_M^c\epsilon$ we get
\begin{alignat}{2} \label{Gcong}
\hat{\G}(\pi_M(1+\pi_M^c\epsilon))&\equiv\pi_L
&&\pmod{\pi_M^m\cdot I_{\G}'},
\end{alignat}
where $I_{\G}'$ is the ideal in $\OO_M[[\epsilon]]$
defined by
\begin{align*}
I_{\G}'&=(\pi_M^{\tilde{\phi}_{M/L}^0(c)}\epsilon^{p^0},
\pi_M^{\tilde{\phi}_{M/L}^1(c)}\epsilon^{p^1},\dots,
\pi_M^{\tilde{\phi}_{M/L}^{\mu}(c)}\epsilon^{p^{\mu}}).
\end{align*}
It follows from (\ref{Fcong}) and (\ref{Gcong})
that there are $r_j,s_k\in R$,
$\delta_{\F}\in(\pi_L,\epsilon)\cdot I_{\F}$, and
$\delta_{\G}\in(\pi_M,\epsilon)\cdot I_{\G}'$ such that
\begin{align} \label{Fu}
\hat{\F}(\pi_L(1+\epsilon))&=\pi_K\cdot
\left(1+\sum_{j=0}^{\nu}\,r_j\pi_L^{i_j}\epsilon^{p^j}
+\delta_{\F}\right) \\
\hat{\G}(\pi_M(1+\pi_M^c\epsilon))&=\pi_L\cdot
\left(1+\sum_{k=0}^{\mu}\,
s_k\pi_M^{\tilde{\phi}_{M/L}^k(c)}\epsilon^{p^k}
+\delta_{\G}\right). \label{Gt}
\end{align}

     Define an ideal in $\OO_M[[\epsilon]]$ by
\begin{align*}
I_{\F\G}&=
\left(\pi_M^{\tilde{\phi}_{L/K}^{j,m}(\tilde{\phi}_{M/L}^k(c))}\epsilon^{p^{j+k}}
:0\le j\le\nu,\;0\le k\le\mu\right) \\
&=\left(\pi_M^{\lambda_{M/K}^g(c)}\epsilon^{p^g}:
0\le g\le\nu+\mu\right).
\end{align*}
Hence for $d\ge0$ and $0\le g\le\nu+\mu$ we have
$\pi_M^d\epsilon^{p^g}\in I_{\F\G}$ if and only if
$d\ge\lambda_{M/K}^g(c)$.  We also define
$u=\pi_L/\pi_M^m\in\OO_M^{\times}$.

\begin{lemma} \label{epspower}
Let $0\le j\le\nu$.  Then
\begin{equation} \nonumber
\pi_L^{i_j}\left(\sum_{k=0}^{\mu}
\,s_k\pi_M^{\tilde{\phi}_{M/L}^k(c)}\epsilon^{p^k}
+\delta_{\G}\right)^{p^j}\equiv u^{i_j}\sum_{k=0}^{\mu}\,
s_k^{p^j}\pi_M^{\tilde{\phi}_{L/K}^{j,m}(\tilde{\phi}_{M/L}^k(c))}
\epsilon^{p^{j+k}}
\pmod{(\pi_M,\epsilon)\cdot I_{\F\G}}.
\end{equation}
\end{lemma}

\proof For $0\le j\le\nu$ define ideals in
$\Z[X_0,X_1,\dots,X_{\mu}]$ by
\begin{align*}
H_j&=(p^hX_k^{p^{j-h}}:1\le h\le j,\;0\le k\le\mu).
\end{align*}
By induction on $j$ we get
\begin{equation} \nonumber
(X_0+X_1+\dots+X_{\mu})^{p^j}\equiv
X_0^{p^j}+X_1^{p^j}+\dots+X_{\mu}^{p^j}
\pmod{H_j}.
\end{equation}
Since both sides of this congruence are homogeneous
polynomials of degree $p^j$, it follows that
\begin{equation} \label{pl}
(X_0+X_1+\dots+X_{\mu})^{p^j}\equiv
X_0^{p^j}+X_1^{p^j}+\dots+X_{\mu}^{p^j}\pmod{H_j'},
\end{equation}
where
\begin{equation} \nonumber
H_j'=(p^hX_k^{p^{j-h}}X_w:1\le h\le j,\;
0\le k\le\mu,\;0\le w\le\mu).
\end{equation}

     Since $\delta_{\G}\in(\pi_M,\epsilon)\cdot I_{\G}'$
there are $\tilde{s}_k\in\OO_M[[\epsilon]]$ such that
$\tilde{s}_k\equiv s_k\pmod{(\pi_M,\epsilon)}$ and
\begin{equation} \nonumber
\sum_{k=0}^{\mu}\,s_k\pi_M^{\tilde{\phi}_{M/L}^k(c)}
\epsilon^{p^k}+\delta_{\G}=\sum_{k=0}^{\mu}\,
\tilde{s}_k\pi_M^{\tilde{\phi}_{M/L}^k(c)}\epsilon^{p^k}.
\end{equation}
Hence by replacing $X_k$ with
$\tilde{s}_k\pi_M^{\tilde{\phi}_{M/L}^k(c)}\epsilon^{p^k}$
for $0\le k\le\mu$ in (\ref{pl}) we get
\begin{equation} \nonumber
\left(\sum_{k=0}^{\mu}\,s_k\pi_M^{\tilde{\phi}_{M/L}^k(c)}
\epsilon^{p^k}+\delta_{\G}\right)^{p^j}\equiv
\sum_{k=0}^{\mu}\,\tilde{s}_k^{p^j}
\pi_M^{p^j\tilde{\phi}_{M/L}^k(c)}\epsilon^{p^{j+k}}
\pmod{\epsilon\cdot A},
\end{equation}
where $A$ is the ideal in $\OO_M[[\epsilon]]$ defined by
\begin{equation} \nonumber
A=(p^h(\pi_M^{\tilde{\phi}_{M/L}^k(c)}
\epsilon^{p^k})^{p^{j-h}}:1\le h\le j,\;0\le k\le\mu).
\end{equation}
Let $1\le h\le j$ and $0\le k\le\mu$.  Since
$i_j+hv_L(p)\ge i_{j-h}$ we have
\begin{align*}
v_M(\pi_L^{i_j}\cdot p^h\pi_M^{p^{j-h}\tilde{\phi}_{M/L}^k(c)})
&\ge mi_{j-h}+p^{j-h}\tilde{\phi}_{M/L}^k(c) \\
&=\tilde{\phi}_{L/K}^{j-h,m}(\tilde{\phi}_{M/L}^k(c)) \\
&\ge\lambda_{M/K}^{j-h+k}(c).
\end{align*}
It follows that $\pi_L^{i_j}\epsilon\cdot
p^h(\pi_M^{\tilde{\phi}_{M/L}^k(c)}\epsilon^{p^k})^{p^{j-h}}
\in\epsilon\cdot I_{\F\G}$, and hence that
$\pi_L^{i_j}\epsilon\cdot A\subset\epsilon\cdot
I_{\F\G}$.  Therefore
\begin{alignat*}{2}
\pi_L^{i_j}\left(\sum_{k=0}^{\mu}
\,s_k\pi_M^{\tilde{\phi}_{M/L}^k(c)}\epsilon^{p^k}
+\delta_{\G}\right)^{p^j}&\equiv
\pi_L^{i_j}\sum_{k=0}^{\mu}\,\tilde{s}_k^{p^j}
\pi_M^{p^j\tilde{\phi}_{M/L}^k(c)}\epsilon^{p^{j+k}}
&&\pmod{\epsilon\cdot I_{\F\G}} \\
&\equiv u^{i_j}\sum_{k=0}^{\mu}\,\tilde{s}_k^{p^j}
\pi_M^{\tilde{\phi}_{L/K}^{j,m}(\tilde{\phi}_{M/L}^k(c))}
\epsilon^{p^{j+k}}&&\pmod{\epsilon\cdot I_{\F\G}}.
\end{alignat*}
Since $\tilde{s}_k\equiv
s_k\pmod{(\pi_M,\epsilon)}$ the lemma follows. \qed
\medskip

     We now replace $\epsilon$ with $\dst\sum_{k=0}^{\mu}
\,s_k\pi_M^{\tilde{\phi}_{M/L}^k(c)}\epsilon^{p^k}
+\delta_{\G}$ in (\ref{Fu}).  With the help of
Lemma~\ref{epspower} we get
\begin{alignat}{2}
\!\hat{\F}(\hat{\G}(\pi_M(1+\pi_M^c\epsilon)))
&=\pi_K\cdot\left(1+\sum_{j=0}^{\nu}
\,r_ju^{i_j}\sum_{k=0}^{\mu}s_k^{p^j}
\pi_M^{\tilde{\phi}_{L/K}^{j,m}(\tilde{\phi}_{M/L}^k(c))}
\epsilon^{p^{j+k}}+\delta_{\F\G}\right) \nonumber \\
&=\pi_K\cdot
\left(1+\sum_{g=0}^{\nu+\mu}
\left(\sum_{(j,k)\in \Omega_g}u^{i_j}r_js_k^{p^j}
\pi_M^{\tilde{\phi}_{L/K}^{j,m}(\tilde{\phi}_{M/L}^k(c))}
\right)\epsilon^{p^g}+\delta_{\F\G}\right)
\label{FG}
\end{alignat}
for some $\delta_{\F\G}\in(\pi_M,\epsilon)\cdot I_{\F\G}$.
\medskip

     To prove (a) in the case $x=c\in\N_0$ we define an
ideal
$J_l=(\pi_M^{nm+\lambda_{M/K}^l(c)},\epsilon^{p^{l+1}})$
in $\OO_M[[\epsilon]]$.  Since
$\pi_K\cdot I_{\F\G}\subset J_l$, by (\ref{FG}) we get
\begin{equation} \nonumber
\hat{\F}(\hat{\G}(\pi_M(1+\pi_M^c\epsilon)))\equiv
\pi_K\pmod{J_l}.
\end{equation}
It follows from Corollary~\ref{ged} that
$\phi_{M/K}^l(c)\ge\lambda_{M/K}^l(c)$.

     Now let $e\ge1$ be relatively prime to
$p[M:K]=pnm$.  Let $\pi_M$ be a uniformizer for $M$, and
choose uniformizers $\pi_L$, $\pi_K$ for $L$, $K$ such
that $\pi_L\equiv\pi_M^m\pmod{\pi_M^{m+1}}$ and
$\pi_K\equiv\pi_L^n\pmod{\pi_L^{n+1}}$; then
$\pi_K\equiv\pi_M^{nm}\pmod{\pi_M^{nm+1}}$.  Let
$\pi_{K_e}\in K^{sep}$ be a root of $X^e-\pi_K$ and set
$K_e=K(\pi_{K_e})$, $L_e=LK_e$, and $M_e=MK_e$.  Let
$0\le h\le\nu$, $0\le i\le\mu$, and $0\le l\le\nu+\mu$.
Then by Lemma~\ref{tame} we get
\begin{align} \label{form1}
\tilde{\phi}_{M/L}^i(x)
&=e^{-1}\tilde{\phi}_{M_e/L_e}^i(ex) \\
\tilde{\phi}_{L/K}^{h,m}(x)
&=e^{-1}\tilde{\phi}_{L_e/K_e}^{h,m}(ex) \\
\phi_{M/K}^l(x)&=e^{-1}\phi_{M_e/K_e}^l(ex)
\label{form3} \\
\lambda_{M/K}^l(x)&=e^{-1}\lambda_{M_e/K_e}^l(ex).
\label{form4}
\end{align}
We know from the preceding paragraph that
$\phi_{M_e/K_e}^l(c)\ge\lambda_{M_e/K_e}^l(c)$ for every
$c\in\N_0$.  By applying (\ref{form3}) and (\ref{form4})
with $x=c/e$
we get $\phi_{M/K}^l(c/e)\ge\lambda_{M/K}^l(c/e)$.
It follows that (a) holds whenever $x=c/e$ with
$c\ge0$, $e\ge1$, and $\gcd(e,pnm)=1$.  Since numbers of
this form are dense in $[0,\infty)$, by continuity we
get $\phi_{M/K}^l(x)\ge\lambda_{M/K}^l(x)$ for all
$x\ge0$.  This proves (a). \medskip

     To facilitate the proof of (b) we define a subset
of the nonnegative reals by
\begin{equation} \label{Tl}
T_l(M/K)=\{t\ge0:\exists\;l_0\le l
\text{ with }|S_l^{l_0}(t)|=1\text{ and }
|S_l^a(t)|=0\text{ for }0\le a< l_0\}.
\end{equation}
Suppose $t>0$ and $(t,\lambda_{M/K}^l(t))$ is not a
vertex of the graph of $\lambda_{M/K}^l$.  Then there is
a unique $0\le l_0\le l$ such that $|S_l^{l_0}(t)|\ge1$;
in fact, $l_0$ is determined by the condition
$(\lambda_{M/K}^l)'(t)=p^{l_0}$.
Hence if the hypotheses of (b) are satisfied with
$x=t$ then $t\in T_l(M/K)$.

\begin{lemma} \label{c0}
Suppose the hypotheses of (b) are satisfied with $x=0$.
Then $0\in T_l(M/K)$.
\end{lemma}

\proof Suppose $0\not\in T_l(M/K)$, and let $l_0$ be the
minimum integer satisfying the hypotheses of (b) with
$x=0$.  Also let $l_1<l_0$ be maximum such that
$|S_l^{l_1}(0)|\not=0$.  Then $|S_l^{l_1}(0)|\ge2$.
Hence there is $(j,k)\in S_l^{l_1}(0)$ such that
$k<\mu$.  Since
\begin{equation} \nonumber
\tilde{\phi}_{M/L}^{k+1}(0)=i_{k+1}'
\le i_k'=\tilde{\phi}_{M/L}^k(0)
\end{equation}
we get
\begin{equation} \nonumber
\lambda_{M/K}^l(0)\le
\tilde{\phi}_{L/K}^{j,m}(\tilde{\phi}_{M/L}^{k+1}(0))\le
\tilde{\phi}_{L/K}^{j,m}(\tilde{\phi}_{M/L}^k(0))
=\lambda_{M/K}^l(0).
\end{equation}
It follows that
$\tilde{\phi}_{L/K}^{j,m}(\tilde{\phi}_{M/L}^{k+1}(0))=
\tilde{\phi}_{L/K}^{j,m}(\tilde{\phi}_{M/L}^k(0))$,
so we have $i_k'=i_{k+1}'$ and
$(j,k+1)\in S_l^{l_1+1}(0)$.  Hence by the maximality of
$l_1$ we get $l_1=l_0-1$.  Since $|S_l^{l_0}(0)|=1$ we
must have $|S_l^{l_0-1}(0)|=2$ and
$(l_0-\mu-1,\mu)\in S_l^{l_0-1}(0)$.  Since
$\tilde{\phi}_{M/L}^{\mu}(0)=0$ we have
\begin{equation} \nonumber
mi_{l_0-\mu}\le mi_{l_0-\mu-1}=
\tilde{\phi}_{L/K}^{l_0-\mu-1,m}(\tilde{\phi}_{M/L}^{\mu}(0))
=\lambda_{M/K}^l(0)\le
\tilde{\phi}_{L/K}^{l_0-\mu,m}(\tilde{\phi}_{M/L}^{\mu}(0))
=mi_{l_0-\mu}
\end{equation}
and hence $\lambda_{M/K}^l(0)=
\tilde{\phi}_{L/K}^{l_0-\mu,m}(\tilde{\phi}_{M/L}^{\mu}(0))$.
Thus $(l_0-\mu,\mu)\in S_l^{l_0}(0)$.
Since $(j,k+1)\in S_l^{l_0}(0)$, and $|S_l^{l_0}(0)|=1$,
we get $k+1=\mu$, and hence
$i_{\mu-1}'=i_k'=i_{k+1}'=i_{\mu}'=0$.  Since
$i_{\mu-1}'>i_{\mu}'=0$, this is a contradiction.  Therefore
$0\in T_l(M/K)$.~\qed

\begin{lemma} \label{rjsk}
Let $c\in\N_0\cap T_l(M/K)$, let $l_0$ be the integer
specified by (\ref{Tl}) for $t=c$, and let $(j,k)$ be
the unique element of $\Omega_{l_0}$ such that
$\lambda_{M/K}^l(c)=
\tilde{\phi}_{L/K}^{j,m}(\tilde{\phi}_{M/L}^k(c))$.
Then $r_j$ and $s_k$ are nonzero.
\end{lemma}

\proof Since $c\in T_l(M/K)$, for $0\le j'<j$ we have
$\tilde{\phi}_{L/K}^{j',m}(\tilde{\phi}_{M/L}^k(c))>
\tilde{\phi}_{L/K}^{j,m}(\tilde{\phi}_{M/L}^k(c))$.  It
follows that $i_{j'}>i_j$, and hence that
$\pi_K\cdot(\pi_L,\epsilon)\cdot
I_{\F}\subset(\pi_L^{n+i_j+1},\epsilon^{p^j+1})$.
Therefore by (\ref{Fu}) we get
\begin{equation} \nonumber
\hat{\F}(\pi_L(1+\epsilon))\equiv\pi_K\cdot
(1+r_j\pi_L^{i_j}\epsilon^{p^j})
\pmod{(\pi_L^{n+i_j+1},\epsilon^{p^j+1})}.
\end{equation}
If $r_j=0$ then by Corollary~\ref{ged} we have
$i_j=\phi_{L/K}^j(0)\ge i_j+1$, a contradiction.  It
follows that $r_j\not=0$.

     Suppose there is $0\le k'<k$ such that
$\tilde{\phi}_{M/L}^{k'}(c)\le\tilde{\phi}_{M/L}^k(c)$.
Since $c\in T_l(M/K)$ we have $(j,k')\not\in
S_l^{j+k'}(c)$, and hence
\begin{equation} \nonumber
\lambda_{M/K}^l(c)<
\tilde{\phi}_{L/K}^j(\tilde{\phi}_{M/L}^{k'}(c))\le
\tilde{\phi}_{L/K}^j(\tilde{\phi}_{M/L}^k(c))
=\lambda_{M/K}^l(c).
\end{equation}
This is a contradiction, so we must have
$\tilde{\phi}_{M/L}^{k'}(c)>\tilde{\phi}_{M/L}^k(c)$ for
$0\le k'<k$.  Hence
$\phi_{M/L}^k(c)=\tilde{\phi}_{M/L}^k(c)$.  Set
$d=\phi_{M/L}^k(c)$.  Then
$\pi_L\cdot(\pi_M,\epsilon)\cdot I_{\G}'\subset
(\pi_M^{m+d+1},\epsilon^{p^k+1})$.  Using (\ref{Gt}) we
get
\begin{equation} \nonumber
\hat{\G}(\pi_M(1+\pi_M^c\epsilon))\equiv\hat{\G}(\pi_M)
(1+s_k\pi_M^d\epsilon^{p^k})
\pmod{(\pi_M^{m+d+1},\epsilon^{p^k+1})}.
\end{equation}
If $s_k=0$ then by Corollary~\ref{ged} we have
$\phi_{M/L}^k(c)\ge d+1$, a contradiction.  It follows
that $s_k\not=0$.~\qed \medskip

     We now prove (b) for $x=c\in\N_0\cap T_l(M/K)$.
Let $l_0$ be the minimum integer satisfying the
hypotheses of (b) for $x=c$.  Then there is a unique
pair $(j,k)\in\Omega_{l_0}$ such that
$\lambda_{M/K}^l(c)=
\tilde{\phi}_{L/K}^{j,m}(\tilde{\phi}_{M/L}^k(c))$.
Furthermore, we have
$\lambda_{M/K}^{l_0}(c)=\lambda_{M/K}^l(c)$ and
$\lambda_{M/K}^{l_1}(c)>\lambda_{M/K}^l(c)$ for
$l_1<l_0$.  Define $J_{l_0}'=
(\pi_M^{nm+\lambda_{M/K}^l(c)+1},\epsilon^{p^{l_0}+1})$.
Then $\pi_L^n\cdot(\pi_M,\epsilon)\cdot
I_{\F\G}\subset J_{l_0}'$, so by (\ref{FG}) we get
\begin{align*}
\hat{\F}(\hat{\G}(\pi_M(1+\pi_M^c\epsilon)))&\equiv
\pi_K\cdot
(1+u^{i_j}r_js_k^{p^j}\pi_M^{\lambda_{M/K}^l(c)}
\epsilon^{p^{l_0}})\pmod{J_{l_0}'}.
\end{align*}
It follows from Lemma~\ref{rjsk} that
$r_j,s_k\in R\smallsetminus\{0\}$ are units.  Therefore
we have
\begin{equation} \nonumber
\hat{\F}(\hat{\G}(\pi_M(1+\pi_M^c\epsilon)))
\not\equiv\pi_K\pmod{J_{l_0}'}.
\end{equation}
Hence by (a) and Corollary~\ref{ged} we get
\begin{equation} \nonumber
\lambda_{M/K}^{l_0}(c)
\le\phi_{M/K}^{l_0}(c)<\lambda_{M/K}^l(c)+1
=\lambda_{M/K}^{l_0}(c)+1.
\end{equation}
Since $\lambda_{M/K}^{l_0}(c)$ and $\phi_{M/K}^{l_0}(c)$ are
integers this implies that
$\lambda_{M/K}^{l_0}(c)=\phi_{M/K}^{l_0}(c)$.  Using (a)
we get
\begin{equation} \nonumber
\lambda_{M/K}^l(c)\le\phi_{M/K}^l(c)\le\phi_{M/K}^{l_0}(c)
=\lambda_{M/K}^{l_0}(c)=\lambda_{M/K}^l(c),
\end{equation}
and hence $\lambda_{M/K}^l(c)=\phi_{M/K}^l(c)$.  Thus
(b) holds for $x\in\N_0\cap T_l(M/K)$.  In particular,
it follows from Lemma~\ref{c0} that (b) holds for $x=0$.

     As in the proof of (a) let $e\ge1$ be relatively
prime to $pnm$, let $\pi_M$ be a uniformizer for $M$,
and choose uniformizers $\pi_L$, $\pi_K$ for $L$, $K$
such that $\pi_L\equiv\pi_M^m\pmod{\pi_M^{m+1}}$ and
$\pi_K\equiv\pi_L^n\pmod{\pi_L^{n+1}}$.  Let
$\pi_{K_e}\in K^{sep}$ be a root of $X^e-\pi_K$ and set
$K_e=K(\pi_{K_e})$, $L_e=LK_e$, and $M_e=MK_e$.
Let $c\in\N_0$ be such that
$c/e\in T_l(M/K)$ and the hypotheses of (b) are
satisfied for the extensions $M/L/K$ with $x=c/e$.
Then it follows from (\ref{form1})--(\ref{form4}) that
$c\in T_l(M_e/K_e)$ and the hypotheses of (b) are
satisfied for the extensions $M_e/L_e/K_e$ with $x=c$.
Hence by the preceding paragraph
we get $\phi_{M_e/K_e}^l(c)=\lambda_{M_e/K_e}^l(c)$.
Using (\ref{form3}) and (\ref{form4}) we deduce that
$\phi_{M/K}^l(c/e)=\lambda_{M/K}^l(c/e)$.

     Now let $r$ be any positive real number such that
the hypotheses of (b) are satisfied with $x=r$, and let
$l_0$ be the minimum integer which satisfies the
hypotheses.  Then there
is a unique element $(j,k)\in \Omega_{l_0}$ such that
$\tilde{\phi}_{L/K}^{j,m}\circ\tilde{\phi}_{M/L}^k(r)
=\lambda_{M/K}^l(r)$.  Let $0\le a\le l_0$ and let
$(u,v)\in \Omega_a$.  Then the graph of
$\tilde{\phi}_{L/K}^{u,m}\circ\tilde{\phi}_{M/L}^v$
is a line of slope $p^{u+v}=p^a\le p^{l_0}$.  Hence if
$(u,v)\not=(j,k)$ and $0\le t<r$ then
$\tilde{\phi}_{L/K}^{u,m}\circ\tilde{\phi}_{M/L}^v(t)>
\tilde{\phi}_{L/K}^{j,m}\circ\tilde{\phi}_{M/L}^k(t)$.
It follows that $S_{l_0}^{l_0}(t)=\{(j,k)\}$ and
$S_{l_0}^a(t)=\varnothing$ for $0\le a<l_0$.  Hence
$t\in T_{l_0}(M/K)$ and the hypotheses of (b) are
satisfied with $x=t$ and $l$ replaced by $l_0$.

     Suppose $\phi_{M/K}^l(r)>\lambda_{M/K}^l(r)$.  Then
there are $c,e\ge1$ such that $\gcd(e,pnm)=1$ and
\begin{equation} \label{bounds}
0<r-\frac{c}{e}<\frac{\phi_{M/K}^l(r)-\lambda_{M/K}^l(r)}{p^{\nu+\mu}}.
\end{equation}
Since $\lambda_{M/K}^{l_0}(r)=\lambda_{M/K}^l(r)$ we get
\begin{equation} \label{g0}
\phi_{M/K}^{l_0}(r)-\lambda_{M/K}^{l_0}(r)\ge
\phi_{M/K}^l(r)-\lambda_{M/K}^l(r)>0.
\end{equation}
Since $\phi_{M/K}^{l_0}$ and $\lambda_{M/K}^{l_0}$ are
continuous increasing piecewise linear functions with
derivatives at most $p^{\nu+\mu}$ it follows from
(\ref{bounds}) and (\ref{g0})
that $\phi_{M/K}^{l_0}(c/e)-\lambda_{M/K}^{l_0}(c/e)>0$.
On the other hand, by the preceding paragraph we know
that $c/e\in T_{l_0}(M/K)$ and the hypotheses of (b) are
satisfied with $x=c/e$ and $l$ replaced by $l_0$.  Hence
$\phi_{M/K}^{l_0}(c/e)=\lambda_{M/K}^{l_0}(c/e)$.  This
is a contradiction, so we must have
$\phi_{M/K}^l(r)\le\lambda_{M/K}^l(r)$.  By combining
this inequality with (a) we get
$\phi_{M/K}^l(r)=\lambda_{M/K}^l(r)$.  This completes
the proof of (b). \medskip

     By setting $x=0$ in Theorem~\ref{ge} we get the
following.  A special case of this result is given in
\cite[Prop.~5.10]{fm}.

\begin{cor}
For $0\le l\le\nu+\mu$ let $i_l''$ denote the $l$th
index of inseparability of $M/K$.  Then
\begin{equation} \nonumber
i_l''\le\min\{mi_j+p^ji_k':(j,k)\in \Omega_{l_0}
\text{ for some }0\le l_0\le l\},
\end{equation}
with equality if there exists $0\le l_0\le l$ such that
there is a unique pair $(j,k)\in \Omega_{l_0}$ which
realizes the minimum.
\end{cor}

\end{document}